\documentclass{article}
\usepackage{ascmac,amsmath,amssymb,amsthm}
\newtheorem{thm}{\bf{Theorem}}

\newtheorem{lem}{\bf{Lemma}}

\begin{document}
\title{The explicit upper bound of the multiple integral of $S(t)$ on the Riemann Hypothesis}
\author{Takahiro Wakasa}
\date{Graduate School of Mathematics, Nagoya University, Chikusa-ku, Nagoya 464-8602, Japan;\\ e-mail: d11003j@math.nagoya-u.ac.jp}
\maketitle
\begin{abstract}
We prove explicit upper bounds of the function $S_m(T)$, defined by the repeated integration of the argument of the Riemann zeta-function. The explicit upper bound of $S(T)$ and $S_1(T)$ have already been obtained by A. Fujii. Our result is a generalization of Fujii's results.
\end{abstract}
\section{Introduction}
We consider the argument of the Riemann zeta function $\zeta(s)$, where $s=\sigma+ti$ is a complex variable, on the critical line $\sigma=\frac{1}{2}$.

We shall give some explicit bounds on $S_m(T)$ defined below under the Riemann hypothesis.

First, we introduce the functions $S(t)$ and $S_1(t)$. When $T$ is not the ordinate of any zero of $\zeta(s)$, we define
\begin{align*}
S(T)=\frac{1}{\pi}\arg\zeta\left(\frac{1}{2}+Ti\right).
\end{align*} 
This is obtained by continuous variation along the straight lines connecting $2$, $2+Ti$, and $\frac{1}{2}+Ti$, starting with the value zero. When $T$ is the ordinate of some zero of $\zeta(s)$, we define
\begin{align*}
S(T)=\frac{1}{2}\{S(T+0)+S(T-0)\}.
\end{align*}

Next, we define $S_1(T)$ by 
\begin{align*}
S_1(T)=\int_{0}^{T}S(t)dt+C,
\end{align*}
where $C$ is the constant defined by
\begin{align*}
C=\frac{1}{\pi}\int_{\frac{1}{2}}^{\infty}\log|\zeta(\sigma)|d\sigma. 
\end{align*}

It is a classical results of von Mangoldt (cf. chapter 9 of Titchmarsh \cite{E.C.TITCHMARSH}) that there exists a number $T_0>0$ such that for $T>T_0$ we have
\begin{align*}
S(T)=O(\log T)
\end{align*}
and
\begin{align*}
S_1(T)=O(\log T).
\end{align*}
Further, it is a classical result of Littlewood \cite{Littlewood} that under the Riemann Hypothesis we  have
\begin{align*}
S(T)=O\left(\frac{\log T}{\log\log T}\right)
\end{align*}
and
\begin{align*}
S_1(T)=O\left(\frac{\log T}{(\log\log T)^2}\right).
\end{align*}

For explicit upper bounds of $|S(T)|$ and $|S_1(T)|$, Karatsuba and Korolev (cf. Theorem 1 and Theorem 2 on \cite{KK}) have shown that
\begin{align*}
|S(T)|<8\log T,
\end{align*}
and
\begin{align*}
|S_1(T)|<1.2\log T
\end{align*}
for $T>T_0$. Also, under the Riemann Hypothesis, it was shown that
\begin{align*}
|S(T)|\leq 0.83\frac{\log T}{\log\log T}
\end{align*} 
for $T>T_0$ in Fujii \cite{Fujii2}. And in Fujii \cite{A. Fujii}, it was shown under the Riemann Hypothesis that
\begin{align*}
|S_1(T)|\leq 0.51\frac{\log T}{(\log\log T)^2}
\end{align*} 
for $T>T_0$.

Next, we introduce the functions $S_2(T)$, $S_3(T)$, $\cdots$ similarly to the case of $S_1(T)$. The non-trivial zeros of $\zeta(s)$ we denote by $\rho=\beta+\gamma i$. When $T\neq\gamma$, we put 
\begin{align*}
S_0(T)=S(T)
\end{align*}
and
\begin{align*}
S_m(T)=\int_{0}^{T}S_{m-1}(t)dt+C_m
\end{align*}
for any integer $m \geq 1$, where $C_m$'s are the constants which are defined by, for any integer $k\geq 1$,
\begin{align*}
C_{2k-1}=\frac{1}{\pi}(-1)^{k-1}\underbrace{\int_{\frac{1}{2}}^{\infty}\int_{\sigma}^{\infty}\cdots \int_{\sigma}^{\infty}}_{(2k-1)-times}\log|\zeta(\sigma)|(d\sigma)^{2k-1},
\end{align*}  
and 
\begin{align*}
C_{2k}=(-1)^{k-1}\underbrace{\int_{\frac{1}{2}}^{\infty}\int_{\sigma}^{\infty}\cdots \int_{\sigma}^{\infty}}_{2k-times}(d\sigma)^{2k}=\frac{(-1)^{k-1}}{(2k)!2^{2k}}.
\end{align*}  
When $T=\gamma$, we put
\begin{align*}
S_m(T)=\frac{1}{2}\{S_m(T+0)+S_m(T-0)\}.
\end{align*}

Concerning $S_m(T)$ for $m \geq 2$, Littlewood \cite{Littlewood} have shown under the Riemann Hypothesis that 
\begin{align*}
S_m(T)=O\left(\frac{\log T}{(\log \log T)^{m+1}}\right).
\end{align*}
 
The purpose of the present article is to prove the following result.
\begin{thm} \label{th}
@\\~~~Under the Riemann Hypothesis for any integer $m \geq 1$, if $m$ is odd,
\begin{align*}
|S_m(t)| &\leq \frac{\log t}{(\log\log t)^{m+1}} \cdot \frac{1}{2\pi m!} \Biggl\{\frac{1}{1-\frac{1}{e}\left(1+\frac{1}{e}\right)}\sum_{j=0}^{m}{}\frac{m!}{(m-j)!}\left(\frac{1}{e}+\frac{1}{2^{j+1}e^2}\right)\\
&~~+\frac{1}{m+1}\cdot \frac{\frac{1}{e}\left(1+\frac{1}{e}\right)}{1-\frac{1}{e}\left(1+\frac{1}{e}\right)}+\frac{1}{m(m+1)}\cdot \frac{1}{1-\frac{1}{e}\left(1+\frac{1}{e}\right)}\Biggr \}\\
&~~+O\left(\frac{\log t}{(\log \log t)^{m+2}}\right).
\end{align*}
~~If $m$ is even,
\begin{align*}
|S_m(t)|&\leq \frac{\log t}{(\log \log t)^{m+1}}\cdot \frac{1}{2\pi m!}\Biggl\{\frac{1}{1-\frac{1}{e}\left(1+\frac{1}{e}\right)}\sum_{j=0}^{m}\frac{m!}{(m-j)!}\left(\frac{1}{e}+\frac{1}{2^{j+1}e^2}\right)\\
&~~+
\frac{1}{m+1}\cdot \frac{\frac{1}{e}\left(1+\frac{1}{e}\right)}{1-\frac{1}{e}\left(1+\frac{1}{e}\right)} +
\frac{\pi}{2}\cdot \frac{1}{1-\frac{1}{e}\left(1+\frac{1}{e}\right)}\Biggr\}+O\left(\frac{\log t}{(\log\log t)^{m+2}}\right).
\end{align*}
\end{thm}
This theorem is a generalization of the known explicit upper bounds for $S(T)$ and $S_1(T)$. It is to be stressed that the argument when the number of integration is odd is different from that when the number of integration is even.

The basic policy of the proof of this theorem is based on A. Fujii \cite{A. Fujii}. In the case when $m$ is odd, we can directly generalize the proof of A. Fujii \cite{A. Fujii}. In the case when $m$ is even, it is an extension of the method of A. Fujii \cite{Fujii2}.

To prove our result, we introduce some more notations. First, we define the function $I_m(T)$ as follows. When $T\neq \gamma$, we put for any integer $k\geq 1$
\begin{align*}
I_{2k-1}(T)=\frac{1}{\pi}(-1)^{k-1}\Re\Biggl\{\underbrace{\int_{\frac{1}{2}}^{\infty}\int_{\sigma}^{\infty}\cdots \int_{\sigma}^{\infty}}_{(2k-1)-times}\log\zeta(\sigma+Ti)(d\sigma)^{2k-1}\Biggr\}
\end{align*}  
and
\begin{align*}
I_{2k}(T)=\frac{1}{\pi}(-1)^{k}\Im\Biggl\{\underbrace{\int_{\frac{1}{2}}^{\infty}\int_{\sigma}^{\infty}\cdots \int_{\sigma}^{\infty}}_{2k-times}\log\zeta(\sigma+Ti)(d\sigma)^{2k}\Biggr\}.
\end{align*}  
When $T=\gamma$, we put for $m\geq 1$
\begin{align*}
I_m(T)=\frac{1}{2}\{I_m(T+0)+I_m(T-0)\}.
\end{align*}
Then, $I_m(T)$ can be expressed as a single integral of the following form (cf. Lemma 2 in Fujii \cite{Fujii3}): for any integer $m\geq 1$
\begin{align*}
I_m(T)=-\frac{1}{\pi}\Im\left\{\frac{i^m}{m!}\int_{\frac{1}{2}}^{\infty}\left(\sigma-\frac{1}{2}\right)^m\frac{\zeta'}{\zeta}(\sigma+Ti)d\sigma \right\}.
\end{align*}
From this expression, it is known under the Riemann Hypothesis that 
\begin{align*}
S_m(T)=I_m(T)
\end{align*}
by Lemma 2 in Fujii \cite{Fujii4}. 

Therefore, we should estimate $I_m(T)$. 

We introduce some lemmas in Section 2 and give the proof of the main results in sections 3 and 4.
\section{Some lemmas}
Here we introduce the following notations.

Let $s=\sigma+ti$. We suppose that $\sigma\geq \frac{1}{2}$ and $t\geq2$. Let $X$ be a positive number satisfying $4\leq X\leq t^2$. Also, we put
\begin{align*}
\sigma_1=\frac{1}{2}+\frac{1}{\log X}
\end{align*}
and
\begin{align*}
\Lambda_X(n)=\left\{ \begin{array}{ll}
\Lambda(n) &~~ {\rm for}~~1\leq n\leq X,  \\
\Lambda(n)\frac{\log\frac{X^2}{n}}{\log X} &~~{\rm for} ~~X\leq n\leq X^2,\\
\end{array} \right.
\end{align*}
with
\begin{align*}
\Lambda(n)=\left\{ \begin{array}{ll}
\log p &~~ {\rm if}~n=p^k~{\rm with~ a~ prime}~p~{\rm and~ an~ integer}~k\geq 1, \\
0&~~{\rm otherwise}.\\
\end{array} \right.
\end{align*}
Using these notations, we state the following lemma.
\begin{lem} \label{lem1}
@\\~~~Let $t\geq 2$, $X>0$ such that $4\leq X \leq t^2$. For $\sigma \geq \sigma_1=\frac{1}{2}+\frac{1}{\log X}$,
\begin{align*}
\frac{\zeta'}{\zeta}(\sigma+ti) &= 
-\sum_{n<X^2}\frac{\Lambda_X(n)}{n^{\sigma+ti}}-
\frac{\left(1+X^{\frac{1}{2}-\sigma} \right) \omega X^{\frac{1}{2}-\sigma}}{1-\frac{1}{e}\left(1+\frac{1}{e}\right)\omega'}\Re\left(\sum_{n<X^2}\frac{\Lambda_X(n)}{n^{\sigma_1+ti}}\right)\\
&~~+\frac{\left(1+X^{\frac{1}{2}-\sigma} \right)\omega X^{\frac{1}{2}-\sigma}}{1-\frac{1}{e}\left(1+\frac{1}{e}\right)\omega'}\cdot\frac{1}{2}\log t
+O\left(X^{\frac{1}{2}-\sigma}\right),
\end{align*}
where $|\omega|\leq 1, -1\leq\omega'\leq 1$.
\end{lem}
This has been proved in Fujii \cite{A. Fujii}. Moreover, we will use the following two lemmas.

\begin{lem} \label{lem2}{\rm (cf. 2.12.7 of Titchmarsh\cite{E.C.TITCHMARSH})}
\begin{align*}
\frac{\zeta'}{\zeta}(s) 
&=\log {2\pi}-1-\frac{E}{2}-\frac{1}{s-1}-\frac{1}{2}\cdot \frac{\Gamma'}{\Gamma}\left(\frac{s}{2}+1\right)+\sum_{\rho}\left(\frac{1}{s-\rho}+\frac{1}{\rho}\right)\\
&=\log {2\pi}-1-\frac{E}{2}-\frac{1}{s-1}-\frac{1}{2} \log \left(\frac{s}{2}+1\right)+\sum_{\rho}\left(\frac{1}{s-\rho}+\frac{1}{\rho}\right)+O\left(\frac{1}{|s|}\right)
\end{align*}
where $E$ is the Euler constant and $\rho$ runs through zeros of $\zeta(s)$.
\end{lem}

\begin{lem} \label{lem3}{\rm (Lemma 1 of Selberg \cite{Selberg2})}
@\\~~~For $X>1$, $s\neq1$, $s\neq -2q~(q=1,2,3,\cdots)$, $s\neq \rho$, 
\begin{align*}
\frac{\zeta'}{\zeta}(s) 
&=-\sum_{n<X^2}\frac{\Lambda_X(n)}{n^s}+\frac{X^{2(1-s)}-X^{1-s}}{(1-s)^2\log X}+\frac{1}{\log X}\sum_{q=1}^{\infty}\frac{X^{-2q-s}-X^{-2(2q+s)}}{(2q+s)^2}\\
&~~~+\frac{1}{\log X}\sum_{\rho}\frac{X^{\rho-s}-X^{2(\rho-s)}}{(s-\rho)^2}.
\end{align*}
\end{lem}
By Lemma \ref{lem2}, we have
\begin{eqnarray}
\Re\frac{\zeta'}{\zeta}(\sigma_1+ti)=-\frac{1}{2}\log t+\sum_{\gamma}\frac{\sigma_1-\frac{1}{2}}{\left(\sigma_1-\frac{1}{2}\right)^2+(t-\gamma)^2}+O(1). \label{use1}
\end{eqnarray}
Since for $\sigma_1\leq \sigma$
\begin{align*}
\frac{1}{\log X}\left|\sum_{\rho}\frac{X^{\rho-s}-X^{2(\rho-s)}}{(s-\rho)^2}\right|&\leq\frac{X^{\frac{1}{2}-\sigma}}{\log X}\sum_{\gamma}\frac{1+X^{\frac{1}{2}-\sigma}}{\left(\sigma-\frac{1}{2}\right)^2+(t-\gamma)^2}\\
&\leq \left(1+X^{\frac{1}{2}-\sigma}\right)X^{\frac{1}{2}-\sigma}\sum_{\gamma}\frac{\sigma_1-\frac{1}{2}}{\left(\sigma_1-\frac{1}{2}\right)^2+(t-\gamma)^2},
\end{align*} 
we have
\begin{align*}
\frac{1}{\log X}\sum_{\rho}\frac{X^{\rho-s}-X^{2(\rho-s)}}{(s-\rho)^2}=\left(1+X^{\frac{1}{2}-\sigma}\right)X^{\frac{1}{2}-\sigma}\cdot \omega \sum_{\gamma}\frac{\sigma_1-\frac{1}{2}}{\left(\sigma_1-\frac{1}{2}\right)^2+(t-\gamma)^2}, 
\end{align*}
where $|\omega|\leq1$. Since for $\sigma\geq \frac{1}{2}$ and $X\leq t^2$
\begin{align*}
\left|\frac{X^{2(1-s)}-X^{1-s}}{(1-s)^2\log X}\right|\ll\frac{X^{2(1-\sigma)}}{t^2\log X}\leq \frac{X^{\frac{1}{2}-\sigma}}{\log X},
\end{align*}
we have for $\sigma_1 \leq \sigma$
\begin{align*}
\frac{\zeta'}{\zeta}(\sigma+ti)&=-\sum_{n<X^2}\frac{\Lambda_X(n)}{n^{\sigma+ti}}+O\left(\frac{X^{\frac{1}{2}-\sigma}}{\log X}\right)\\
&~~~+\left(1+X^{\frac{1}{2}-\sigma}\right)\omega X^{\frac{1}{2}-\sigma}\sum_{\gamma}\frac{\sigma_1-\frac{1}{2}}{\left(\sigma_1-\frac{1}{2}\right)^2+(t-\gamma)^2}
\end{align*}
by Lemma \ref{lem3}. Especially,
\begin{align}
\Re\frac{\zeta'}{\zeta}(\sigma_1+ti)&=\Re\left(\sum_{n<X^2}\frac{\Lambda_X(n)}{n^{\sigma_1+ti}}\right)+O\left(\frac{1}{\log X}\right)\nonumber \\
&~~~+\left(1+\frac{1}{e}\right)\frac{1}{e}\omega'\sum_{\gamma}\frac{\sigma_1-\frac{1}{2}}{\left(\sigma_1-\frac{1}{2}\right)^2+(t-\gamma)^2},\label{use2}
\end{align}
where $-1\leq \omega' \leq 1$.

Hence by (\ref{use1}) and (\ref{use2}), we get
\begin{eqnarray}
\sum_{\gamma}\frac{\sigma_1-\frac{1}{2}}{\left(\sigma_1-\frac{1}{2}\right)^2+(t-\gamma)^2}=\frac{1}{1-\frac{1}{e}\left(1+\frac{1}{e}\right)\omega'}\cdot\frac{1}{2}\log t+O\left(\left|\sum_{n<X^2}\frac{\Lambda_X(n)}{n^{\sigma_1+ti}}\right|\right).\label{1.9}
\end{eqnarray}
This relation will be used in the following proof of Theorem \ref{th}.
\section{Proof of Theorem 1 in the case when $m$ is odd}
If $m$ is odd, we have
\begin{eqnarray}
I_m(t)&=&\frac{i^{m+1}}{\pi m!}\Im \Biggl\{i\Biggl\{\int_{\sigma_1}^{\infty}\left(\sigma-\frac{1}{2}\right)^m  \frac{\zeta'}{\zeta}(\sigma+ti)d\sigma
+\frac{\left(\sigma_1-\frac{1}{2}\right)^{m+1}}{m+1}\cdot \frac{\zeta'}{\zeta}(\sigma_1+ti)\nonumber \\
&~&-\int_{\frac{1}{2}}^{\sigma_1}\left(\sigma-\frac{1}{2}\right)^m\left\{\frac{\zeta'}{\zeta}(\sigma_1+ti)-\frac{\zeta'}{\zeta}(\sigma+ti)\right\}d\sigma \Biggr\}\Biggr\}\nonumber \\
&=&\frac{i^{m+1}}{\pi m!}\Im \left\{i(J_1+J_2+J_3)\right\},\label{1.1}
\end{eqnarray}
say.

First, we estimate $J_1$. By Lemma \ref{lem1},
\begin{align}
J_1 &=\int_{\sigma_1}^{\infty}\left(\sigma-\frac{1}{2}\right)^m \Biggl\{-\sum_{n<X^2}\frac{\Lambda_X(n)}{n^{\sigma+ti}}-\frac{\left(1+X^{\frac{1}{2}-\sigma} \right) \omega X^{\frac{1}{2}-\sigma}}{1-\frac{1}{e}\left(1+\frac{1}{e}\right)\omega'}\Re\left(\sum_{n<X^2}\frac{\Lambda_X(n)}{n^{\sigma_1+ti}}\right)\nonumber \\
&~~+\frac{\left(1+X^{\frac{1}{2}-\sigma} \right)\omega X^{\frac{1}{2}-\sigma}}{1-\frac{1}{e}\left(1+\frac{1}{e}\right)\omega'}\cdot\frac{1}{2}\log t+O\left(X^{\frac{1}{2}-\sigma}\right)\Biggr\}d\sigma \nonumber \\
&=-\int_{\sigma_1}^{\infty}\left(\sigma-\frac{1}{2}\right)^m \sum_{n<X^2}\frac{\Lambda_X(n)}{n^{\sigma+ti}}d\sigma+\eta_1(t),\nonumber
\end{align}
say. Then, by integration by parts repeatedly 
\begin{align}
J_1=-\sum_{j=0}^{m}\left(\frac{m!}{(m-j)!}\left(\sigma_1-\frac{1}{2}\right)^{m-j}\sum_{n<X^2}\frac{\Lambda_X(n)}{n^{\sigma_1+ti}(\log n)^{j+1}}\right)+\eta_1(t).\label{1.2}
\end{align} 
And we have
\begin{align}
|\eta_1(t)| &=\left|\int_{\sigma_1}^{\infty}\left(\sigma-\frac{1}{2}\right)^m\frac{\left(1+X^{\frac{1}{2}-\sigma} \right) \omega X^{\frac{1}{2}-\sigma}}{1-\frac{1}{e}\left(1+\frac{1}{e}\right)\omega'}d\sigma \right| \cdot \left|-\Re\left(\sum_{n<X^2}\frac{\Lambda_X(n)}{n^{\sigma_1+ti}}\right)
+\frac{1}{2}\log t\right|\nonumber \\
&~~+O\left\{\int_{\sigma_1}^{\infty}\left(\sigma-\frac{1}{2}\right)^m X^{\frac{1}{2}-\sigma}d\sigma \right\}\nonumber \\
&\leq \frac{1}{1-\frac{1}{e}\left(1+\frac{1}{e}\right)}\left|\frac{1}{2}\log t-\Re\left(\sum_{n<X^2}\frac{\Lambda_X(n)}{n^{\sigma_1+ti}}\right)
\right| \nonumber \\
&~~\cdot \int_{\sigma_1}^{\infty}\left(\sigma-\frac{1}{2}\right)^m \left(1+X^{\frac{1}{2}-\sigma}\right)X^{\frac{1}{2}-\sigma}d\sigma+O\left\{\int_{\sigma_1}^{\infty}\left(\sigma-\frac{1}{2}\right)^m X^{\frac{1}{2}-\sigma}d\sigma \right\}\nonumber \\
&\leq \frac{1}{1-\frac{1}{e}\left(1+\frac{1}{e}\right)}\cdot \frac{1}{2}\log t\cdot \frac{1}{(\log X)^{m+1}}\left(\sum_{j=0}^{m}\frac{m!}{(m-j)!}\left(\frac{1}{e}+\frac{1}{2^{j+1}e^2}\right)\right)\nonumber \\
&~~+O\left(\frac{1}{(\log X)^{m+1}}\left|\sum_{n<X^2}\frac{\Lambda_X(n)}{n^{\sigma_1+ti}}\right|\right)\nonumber \\
&=\eta_2(t)+O\left(\frac{1}{(\log X)^{m+1}}\left|\sum_{n<X^2}\frac{\Lambda_X(n)}{n^{\sigma_1+ti}}\right|\right), \label{1.3}
\end{align}
say, since by partial integration
\begin{align*}
\int_{\sigma_1}^{\infty}\left(\sigma-\frac{1}{2}\right)^m \left(1+X^{\frac{1}{2}-\sigma}\right)X^{\frac{1}{2}-\sigma}d\sigma =
\frac{1}{(\log X)^{m+1}}\left(\sum_{j=0}^{m}\frac{m!}{(m-j)!}\left(\frac{1}{e}+\frac{1}{2^{j+1}e^2}\right)\right).
\end{align*}

Next, applying Lemma \ref{lem1} to $J_2$, we get
\begin{align}
J_2&=\frac{1}{(m+1)(\log X)^{m+1}}\cdot\Biggl\{\sum_{n<X^2}\frac{\Lambda_X(n)}{n^{\sigma_1+ti}}-\frac{\left(1+\frac{1}{e}\right) \frac{1}{e}\omega}{1-\frac{1}{e}\left(1+\frac{1}{e}\right)\omega' }\Re \left(\sum_{n<X^2}\frac{\Lambda_X(n)}{n^{\sigma_1+ti}}\right)\nonumber \\
&~~+\frac{\left(1+\frac{1}{e}\right) \frac{1}{e}\omega}{1-\frac{1}{e}\left(1+\frac{1}{e}\right)\omega' }\cdot \frac{1}{2}\log t+O\left(X^{\frac{1}{2}-\sigma_1}\right)\Biggr\}\nonumber \\
&=\frac{1}{(m+1)(\log X)^{m+1}}\cdot \frac{\left(1+\frac{1}{e}\right)\frac{1}{e}\omega}{1-\frac{1}{e}\left(1+\frac{1}{e}\right)\omega' }\cdot \frac{1}{2}\log t+O\left\{\frac{1}{(\log X)^{m+1}}\left|\sum_{n<X^2}\frac{\Lambda_X(n)}{n^{\sigma_1+ti}}\right|\right\}\nonumber \\
&=\eta_3(t)+O\left\{\frac{1}{(\log X)^{m+1}}\left|\sum_{n<X^2}\frac{\Lambda_X(n)}{n^{\sigma_1+ti}}\right|\right\}, \label{1.4}
\end{align}
say.

Next, we estimate $J_3$. By Lemma \ref{lem2}, we have
\begin{align}
\Im(iJ_3) &=\Re(J_3) =-\int_{\frac{1}{2}}^{\sigma_1} \left( \sigma-\frac{1}{2}\right )^m \Re \left\{\frac{\zeta'}{\zeta}(\sigma_1+ti)- \frac{\zeta'}{\zeta}(\sigma+ti) \right\} d\sigma \nonumber \\
&=-\sum_{\gamma}\frac{1}{\left(\sigma_1-\frac{1}{2}\right)^2+(t-\gamma)^2}\cdot \nonumber \\
&~~\int_{\frac{1}{2}}^{\sigma_1}\left(\sigma-\frac{1}{2}\right)^m \frac{(\sigma_1-\sigma) \left\{(t-\gamma )^2-\left(\sigma_1-\frac{1}{2}\right) \left(\sigma-\frac{1}{2}\right)\right\}}{\left(\sigma-\frac{1}{2}\right)^2+(t-\gamma)^2}d\sigma \nonumber \\
&~~+O\left(\frac{1}{t(\log X)^{m+1}}\right) \nonumber \\
&=-\sum_{\gamma}\frac{1}{\left(\sigma_1-\frac{1}{2}\right)^2+(t-\gamma)^2}\cdot K(\gamma)+O\left(\frac{1}{(\log X)^{m+1}}\right), \label{1.5}
\end{align}
say, where $\gamma$ is the imaginary part of $\rho=\beta+\gamma i$.

If $t=\gamma$,
\begin{align}
K(\gamma)&=-\int_{\frac{1}{2}}^{\sigma_1}\left(\sigma-\frac{1}{2}\right)^{m-1}\left(\sigma_1-\frac{1}{2}\right)\left(\sigma_1-\sigma\right)d\sigma \nonumber \\
&=-\frac{1}{m(m+1)}\left(\sigma_1-\frac{1}{2}\right)^{m+2}.\label{1.6}
\end{align}

If $t\neq\gamma$, by putting $\sigma-\frac{1}{2}=v $, $\sigma_1-\frac{1}{2}=\frac{1}{\log X}=\Delta$ and $|t-\gamma|=B$, we get
\begin{align*}
K(\gamma)&=\int_{0}^{\Delta}v^m\frac{(\Delta-v)(B^2-\Delta v)}{v^2+B^2}dv\\
&=\int_{0}^{\Delta}\left\{v^m \Delta -(B^2+\Delta)v^{m-1} +\frac{B^2(B^2+\Delta ^2)v^{m-1}}{v^2+B^2}\right\}dv\\
&=\frac{\Delta^{m+2}}{m+1}-\frac{(B^2+\Delta^2)\Delta^m}{m}+\int_{0}^{\Delta}\frac{(B^2+\Delta^2)v^{m-1}}{\left(\frac{v}{B}\right)^2+1}dv.
\end{align*}
Putting $\frac{v}{B}=u$, we have
\begin{align*}
K(\gamma)&=\frac{\Delta^{m+2}}{m+1}-\frac{(B^2+\Delta^2)\Delta^m}{m}+(B^2+\Delta^2)\int_{0}^{\frac{\Delta}{B}}\frac{(uB)^{m-1}B}{1+u^2}du\\
&=\frac{\Delta^{m+2}}{m+1}-\frac{(B^2+\Delta^2)\Delta^m}{m}\\
&~~+(B^2+\Delta^2)B^m i^{m+1}\left\{\sum_{j=1}^{\frac{m-1}{2}}\frac{(-1)^{j-1}}{2j-1}\left(\frac{\Delta}{B}\right)^{2j-1}-\arctan \left(\frac{\Delta}{B}\right)\right\}\\
&=\Delta^{m+2}\Biggl\{\frac{1}{m+1}-\frac{1}{m}\cdot \frac{B^2}{\Delta^2}-\frac{1}{m}\\
&~~+\left(\frac{B^{m+2}}{\Delta^{m+2}}+\frac{B^m}{\Delta^m}\right) i^{m+1}\left.\left\{\sum_{j=1}^{\frac{m-1}{2}}\frac{(-1)^{j-1}}{2j-1}\left(\frac{\Delta}{B}\right)^{2j-1}-\arctan \left(\frac{\Delta}{B}\right)\right\} \right\}.
\end{align*}
Putting $y=\frac{\Delta}{B}$, we get
\begin{align}
K(\gamma)&=\Delta^{m+2}\Biggl\{-i^{m+1}\left(\frac{1}{y^{m+2}}+\frac{1}{y^m}\right)\arctan y-\frac{1}{my^2}\nonumber \\
&~~+i^{m+1}\left(\frac{1}{y^{m+2}}+\frac{1}{y^m}\right)\sum_{j=1}^{\frac{m-1}{2}}\frac{(-1)^{j-1}}{2j-1}y^{2j-1}-\frac{1}{m(m+1)}\Biggr\}\nonumber\\
&=\Delta^{m+2}\left(g(y)-\frac{1}{m(m+1)}\right).\label{1.7}
\end{align}
When $y$ tends to $0$, $g(y)$ is convergent to $\frac{2}{m(m+2)}$ since
\begin{align*}
g(y)=\frac{2}{m(m+2)}-\frac{2}{(m+2)(m+4)}y^2+\frac{2}{(m+4)(m+6)}y^4-\cdots.
\end{align*}
When $y$ tends to infinity, $g(y)$ tends to $0$. Hence for $y>0$, we get $g'(y)<0$. 
Hence
\begin{eqnarray*}
0\leq g(y)\leq \frac{2}{m(m+2)},
\end{eqnarray*}
so that
\begin{eqnarray}
-\frac{1}{m(m+1)}\leq g(y)-\frac{1}{m(m+1)}\leq \frac{1}{(m+1)(m+2)}.\label{A}
\end{eqnarray}
Therefore by (\ref{1.7}) and (\ref{A}), we obtain
\begin{eqnarray*}
-\frac{\Delta^{m+2}}{m(m+1)}\leq K(\gamma) \leq \frac{\Delta^{m+2}}{(m+1)(m+2)},
\end{eqnarray*}
so that
\begin{eqnarray*}
-\frac{1}{m(m+1)}\left(\sigma_1-\frac{1}{2}\right)^{m+2}\leq K(\gamma) \leq \frac{1}{(m+1)(m+2)}\left(\sigma_1-\frac{1}{2}\right)^{m+2}.
\end{eqnarray*}
Hence
\begin{eqnarray}
-\sum_{\gamma}\frac{1}{\left(\sigma_1-\frac{1}{2}\right)^2+(t-\gamma)^2}K(\gamma)\leq \frac{\left(\sigma_1-\frac{1}{2}\right)^{m+2}}{m(m+1)}\sum_{\gamma}\frac{1}{\left(\sigma_1-\frac{1}{2}\right)^2+(t-\gamma)^2}\label{B}
\end{eqnarray}
and
\begin{eqnarray}
-\sum_{\gamma}\frac{1}{\left(\sigma_1-\frac{1}{2}\right)^2+(t-\gamma)^2}K(\gamma)\geq -\frac{\left(\sigma_1-\frac{1}{2}\right)^{m+2}}{(m+1)(m+2)}\sum_{\gamma}\frac{1}{\left(\sigma_1-\frac{1}{2}\right)^2+(t-\gamma)^2}.\label{C}
\end{eqnarray}
By (\ref{1.9}), (\ref{B}) and (\ref{C}), we have
\begin{align}
-\sum_{\gamma}\frac{1}{\left(\sigma_1-\frac{1}{2}\right)^2+(t-\gamma)^2}K(\gamma)&\leq 
\frac{\left(\sigma_1-\frac{1}{2}\right)^{m+1}}{m(m+1)}\Biggl\{\frac{1}{1-\frac{1}{e}\left(1+\frac{1}{e}\right)\omega'}\cdot\frac{1}{2}\log t\nonumber \\
&~~~+O\left(\left|\sum_{n<X^2}\frac{\Lambda_X(n)}{n^{\sigma_1+ti}}\right|\right)\Biggr\}\label{D}
\end{align}
and
\begin{align}
-\sum_{\gamma}\frac{1}{\left(\sigma_1-\frac{1}{2}\right)^2+(t-\gamma)^2}K(\gamma)&\geq 
-\frac{\left(\sigma_1-\frac{1}{2}\right)^{m+1}}{(m+1)(m+2)}\Biggl\{\frac{1}{1-\frac{1}{e}\left(1+\frac{1}{e}\right)\omega'}\cdot\frac{1}{2}\log t\nonumber \\
&~~~+O\left(\left|\sum_{n<X^2}\frac{\Lambda_X(n)}{n^{\sigma_1+ti}}\right|\right)\Biggr\}.\label{E}
\end{align}
Hence by (\ref{1.5}), (\ref{1.6}), (\ref{D}) and (\ref{E}), if $m\equiv 1 \pmod 4$,
\begin{align}
i^{m+1}\Im(iJ_3)&\leq \frac{1}{(m+1)(m+2)}\cdot\frac{1}{(\log X)^{m+1}}\cdot \frac{1}{1-\frac{1}{e}\left(1+\frac{1}{e}\right)\omega'}\cdot\frac{1}{2}\log t\nonumber \\
&~~~+O\left(\frac{1}{(\log X)^{m+1}}\left|\sum_{n<X^2}\frac{\Lambda_X(n)}{n^{\sigma_1+ti}}\right|\right)\nonumber \\
&=\eta_4(t)+O\left(\frac{1}{(\log X)^{m+1}}\left|\sum_{n<X^2}\frac{\Lambda_X(n)}{n^{\sigma_1+ti}}\right|\right)\label{1.10}
\end{align}
and
\begin{align}
i^{m+1}\Im(iJ_3)&\geq -\frac{1}{m(m+1)}\cdot\frac{1}{(\log X)^{m+1}}\cdot \frac{1}{1-\frac{1}{e}\left(1+\frac{1}{e}\right)\omega'}\cdot\frac{1}{2}\log t\nonumber \\
&~~~+O\left(\frac{1}{(\log X)^{m+1}}\left|\sum_{n<X^2}\frac{\Lambda_X(n)}{n^{\sigma_1+ti}}\right|\right)\nonumber \\
&=-\eta_5(t)+O\left(\frac{1}{(\log X)^{m+1}}\left|\sum_{n<X^2}\frac{\Lambda_X(n)}{n^{\sigma_1+ti}}\right|\right)\label{1.11}
\end{align}
since $-i^{m+1}=1$ and $i^{m+1}=-1$. And if $m\equiv 3 \pmod 4$,
\begin{align}
i^{m+1}\Im(iJ_3)&\leq \frac{1}{m(m+1)}\cdot\frac{1}{(\log X)^{m+1}}\cdot \frac{1}{1-\frac{1}{e}\left(1+\frac{1}{e}\right)\omega'}\cdot\frac{1}{2}\log t\nonumber \\
&~~~+O\left(\frac{1}{(\log X)^{m+1}}\left|\sum_{n<X^2}\frac{\Lambda_X(n)}{n^{\sigma_1+ti}}\right|\right)\nonumber \\
&=\eta_5(t)+O\left(\frac{1}{(\log X)^{m+1}}\left|\sum_{n<X^2}\frac{\Lambda_X(n)}{n^{\sigma_1+ti}}\right|\right)\label{1.12}
\end{align}
and
\begin{align}
i^{m+1}\Im(iJ_3)&\geq -\frac{1}{(m+1)(m+2)}\cdot\frac{1}{(\log X)^{m+1}}\cdot \frac{1}{1-\frac{1}{e}\left(1+\frac{1}{e}\right)\omega'}\cdot\frac{1}{2}\log t\nonumber \\
&~~~+O\left(\frac{1}{(\log X)^{m+1}}\left|\sum_{n<X^2}\frac{\Lambda_X(n)}{n^{\sigma_1+ti}}\right|\right)\nonumber \\
&=-\eta_4(t)+O\left(\frac{1}{(\log X)^{m+1}}\left|\sum_{n<X^2}\frac{\Lambda_X(n)}{n^{\sigma_1+ti}}\right|\right)\label{1.13}
\end{align}
since $-i^{m+1}=-1$ and $i^{m+1}=1$.

Therefore by (\ref{1.1}), (\ref{1.3}), (\ref{1.4}), (\ref{1.10}), (\ref{1.11}), (\ref{1.12}) and (\ref{1.13}), we obtain
\begin{align}
I_m(t)&=\frac{1}{\pi m!}\Biggl\{-i^{m+1}\sum_{j=0}^{m}\left(\frac{m!}{(m-j)!}\left(\sigma_1-\frac{1}{2}\right)^{m-j}\sum_{n<X^2}\frac{\Lambda_X(n)}{n^{\sigma_1+ti}(\log n)^{j+1}}\right)\nonumber \\
&~~~+O\left(\frac{1}{(\log X)^{m+1}}\left|\sum_{n<X^2}\frac{\Lambda_X(n)}{n^{\sigma_1+ti}}\right|\right)\Biggr\}\nonumber \\
&~~~+\frac{1}{\pi m!}\cdot \Xi(t), \label{1.14}
\end{align}
where $\Xi(t)$ satisfies the following inequalities. If $m\equiv 1 \pmod 4$,
\begin{align*}
\Xi(t)&\leq \eta_2(t)-\eta_3(t)+\eta_4(t)\\
&=\frac{1}{1-\frac{1}{e}\left(1+\frac{1}{e}\right)}\cdot\frac{1}{2}\log t\cdot \frac{1}{(\log X)^{m+1}}\left(\sum_{j=0}^{m}\frac{m!}{(m-j)!}\left(\frac{1}{e}+\frac{1}{2^{j+1}e^2}\right)\right)\\
&~~~-\frac{1}{m+1}\cdot \frac{\left(1+\frac{1}{e}\right) \frac{1}{e}\omega}{1-\frac{1}{e}\left(1+\frac{1}{e}\right)\omega'}\cdot\frac{1}{2}\log t\cdot\frac{1}{(\log X)^{m+1}}\\
&~~~+\frac{1}{(m+1)(m+2)}\cdot\frac{1}{1-\frac{1}{e}\left(1+\frac{1}{e}\right)\omega'}\cdot\frac{1}{2}\log t\cdot\frac{1}{(\log X)^{m+1}},
\end{align*}
and
\begin{align*}
\Xi(t)&\geq -\eta_2(t)-\eta_3(t)-\eta_5(t)\\
&=-\frac{1}{1-\frac{1}{e}\left(1+\frac{1}{e}\right)}\cdot\frac{1}{2}\log t\cdot \frac{1}{(\log X)^{m+1}}\left(\sum_{j=0}^{m}\frac{m!}{(m-j)!}\left(\frac{1}{e}+\frac{1}{2^{j+1}e^2}\right)\right)\\
&~~~-\frac{1}{m+1}\cdot \frac{\left(1+\frac{1}{e}\right) \frac{1}{e}\omega}{1-\frac{1}{e}\left(1+\frac{1}{e}\right)\omega'}\cdot\frac{1}{2}\log t\cdot\frac{1}{(\log X)^{m+1}}\\
&~~~-\frac{1}{m(m+1)}\cdot\frac{1}{1-\frac{1}{e}\left(1+\frac{1}{e}\right)\omega'}\cdot\frac{1}{2}\log t\cdot\frac{1}{(\log X)^{m+1}},
\end{align*}
and if $m\equiv 3 \pmod 4$,
\begin{align*}
\Xi(t)&\leq \eta_2(t)+\eta_3(t)+\eta_5(t)\\
&=\frac{1}{1-\frac{1}{e}\left(1+\frac{1}{e}\right)}\cdot\frac{1}{2}\log t\cdot \frac{1}{(\log X)^{m+1}}\left(\sum_{j=0}^{m}\frac{m!}{(m-j)!}\left(\frac{1}{e}+\frac{1}{2^{j+1}e^2}\right)\right)\\
&~~~+\frac{1}{m+1}\cdot \frac{\left(1+\frac{1}{e}\right) \frac{1}{e}\omega}{1-\frac{1}{e}\left(1+\frac{1}{e}\right)\omega'}\cdot\frac{1}{2}\log t\cdot\frac{1}{(\log X)^{m+1}}\\
&~~~+\frac{1}{m(m+1)}\cdot\frac{1}{1-\frac{1}{e}\left(1+\frac{1}{e}\right)\omega'}\cdot\frac{1}{2}\log t\cdot\frac{1}{(\log X)^{m+1}},
\end{align*}
and
\begin{align*}
\Xi(t)&\geq -\eta_2(t)-\eta_3(t)-\eta_4(t)\\
&=-\frac{1}{1-\frac{1}{e}\left(1+\frac{1}{e}\right)}\cdot\frac{1}{2}\log t\cdot \frac{1}{(\log X)^{m+1}}\left(\sum_{j=0}^{m}\frac{m!}{(m-j)!}\left(\frac{1}{e}+\frac{1}{2^{j+1}e^2}\right)\right)\\
&~~~-\frac{1}{m+1}\cdot \frac{\left(1+\frac{1}{e}\right) \frac{1}{e}\omega}{1-\frac{1}{e}\left(1+\frac{1}{e}\right)\omega'}\cdot\frac{1}{2}\log t\cdot\frac{1}{(\log X)^{m+1}}\\
&~~~-\frac{1}{(m+1)(m+2)}\cdot\frac{1}{1-\frac{1}{e}\left(1+\frac{1}{e}\right)\omega'}\cdot\frac{1}{2}\log t\cdot\frac{1}{(\log X)^{m+1}}.
\end{align*}
In (\ref{1.14}), we have
\begin{align}
\left|\sum_{n<X^2}\frac{\Lambda_X(n)}{n^{\sigma_1+ti}}\right|\leq \sum_{n<X}\frac{\Lambda(n)}{n^\frac{1}{2}}+\sum_{X\leq n\leq X^2}\frac{\Lambda(n)\log\frac{X^2}{n}}{n^\frac{1}{2}}\cdot\frac{1}{\log X} \ll \frac{X}{\log X}.\label{sigxx}
\end{align}
Hence the second term on the right-hand side of (\ref{1.14}) is $ \ll \frac{X}{(\log X)^{m+2}}$. Similarly, since
\begin{align*}
\left|\sum_{n<X^2}\frac{\Lambda_X(n)}{n^{\sigma_1+ti}(\log n)^{j+1}}\right|&\leq \sum_{n<X}\frac{\Lambda(n)}{n^\frac{1}{2}(\log n)^{j+1}}+\sum_{X\leq n\leq X^2}\frac{\Lambda(n)\log\frac{X^2}{n}}{n^\frac{1}{2}(\log n)^{j+1}}\cdot\frac{1}{\log X}\\
&\ll \frac{X}{(\log X)^{j+2}},
\end{align*}
we estimate that the first term on the right-hand side of (\ref{1.14}) is $ \ll \frac{X}{(\log X)^{m+2}}$.

Therefore, taking $ X=\log t$, we obtain
\begin{align*}
|I_m(t)|&= \frac{1}{\pi m!}\Xi(t)+O\left(\frac{\log t}{(\log\log t)^{m+2}}\right)\\
&= \frac{\log t}{(\log\log t)^{m+1}} \cdot \frac{1}{2\pi m!} \Biggl\{\frac{1}{1-\frac{1}{e}\left(1+\frac{1}{e}\right)}\sum_{j=0}^{m}\frac{m!}{(m-j)!}\left(\frac{1}{e}+\frac{1}{2^{j+1}e^2}\right)\\
&~~+\frac{1}{m+1}\cdot \frac{\frac{1}{e}\left(1+\frac{1}{e}\right)}{1-\frac{1}{e}\left(1+\frac{1}{e}\right)}+\frac{1}{m(m+1)}\cdot \frac{1}{1-\frac{1}{e}\left(1+\frac{1}{e}\right)}\Biggr \}\\
&~~+O\left(\frac{\log t}{(\log\log t)^{m+2}}\right).
\end{align*}
This is the first part of the theorem.

\section{Proof of Theorem 1 in the case when $m$ is even}
If $m$ is even, we get similarly
\begin{eqnarray}
I_m(t)&=&\frac{-i^{m}}{\pi m!}\Im \Biggl\{\Biggl\{\int_{\sigma_1}^{\infty}\left(\sigma-\frac{1}{2}\right)^m  \frac{\zeta'}{\zeta}(\sigma+ti)d\sigma
+\frac{\left(\sigma_1-\frac{1}{2}\right)^{m+1}}{m+1}\cdot \frac{\zeta'}{\zeta}(\sigma_1+ti)\nonumber \\
&~&-\int_{\frac{1}{2}}^{\sigma_1}\left(\sigma-\frac{1}{2}\right)^m\left\{\frac{\zeta'}{\zeta}(\sigma_1+ti)-\frac{\zeta'}{\zeta}(\sigma+ti)\right\}d\sigma \Biggr\}\Biggr\}\nonumber \\
&=&\frac{-i^{m}}{\pi m!}\Im \left\{(J_1+J_2+J_3)\right\},\label{1.15}
\end{eqnarray}
say.
By Lemma \ref{lem1} and (\ref{sigxx}), we have 
\begin{align}
J_1&=-\int_{\sigma_1}^{\infty}\left(\sigma-\frac{1}{2}\right)^m\sum_{n<X^2}\frac{\Lambda_X(n)}{n^{\sigma+ti}}d\sigma+\int_{\sigma_1}^{\infty}\left(\sigma-\frac{1}{2}\right)^mO\left(X^{\frac{1}{2}-\sigma}\right)d\sigma \nonumber\\
&~~~~~~~+\int_{\sigma_1}^{\infty}\left(\sigma-\frac{1}{2}\right)^m\Biggl\{-\frac{\left(1+X^{\frac{1}{2}-\sigma}\right)\omega X^{\frac{1}{2}-\sigma}}{1-\frac{1}{e}\left(1+\frac{1}{e}\right)\omega'}\Re\left(\sum_{n<X^2}\frac{\Lambda_X(n)}{n^{\sigma_1+ti}}\right)\nonumber \\
&~~~~~~~~~~~~~~~~~~~~~~~~~~~~~~~~~~~~~+\frac{\left(1+X^{\frac{1}{2}-\sigma}\right)\omega X^{\frac{1}{2}-\sigma}}{1-\frac{1}{e}\left(1+\frac{1}{e}\right)\omega'}\cdot\frac{1}{2}\log t \Biggr\} d\sigma \nonumber\\
&=-\int_{\sigma_1}^{\infty}\left(\sigma-\frac{1}{2}\right)^m\sum_{n<X^2}\frac{\Lambda_X(n)}{n^{\sigma+ti}}d\sigma +O\left\{\int_{\sigma_1}^{\infty}\left(\sigma-\frac{1}{2}\right)^mX^{\frac{1}{2}-\sigma}d\sigma\right\}+\eta'_1(t)\nonumber\\
&=\sum_{j=0}^{m}\frac{m!}{(m-j)!}\left(\sigma_1+\frac{1}{2}\right)^{m-j}\sum_{n<X^2}\frac{\Lambda_X(n)}{n^{\sigma_1+ti}(\log n)^{j+1}}\nonumber\\
&~~~~~~~~~~~~~~~~~~~~~~~~~~~~~~~~~~~+O\left(\frac{1}{(\log X)^{m+1}}\left|\sum_{n<X^2}\frac{\Lambda_X(n)}{n^{\sigma_1+ti}}\right|\right)+\eta'_1(t)\nonumber \\
&\ll \frac{X}{(\log X)^{m+2}}+\eta'_1(t),\label{1.16}
\end{align}
say, and
\begin{align}
J_2&=\frac{1}{(m+1)(\log X)^{m+1}}\Biggl\{\sum_{n<X^2}\frac{\Lambda_X(n)}{n^{\sigma_1+ti}}-\frac{\left(1+\frac{1}{e}\right)\frac{1}{e}\omega}{1-\frac{1}{e}\left(1+\frac{1}{e}\right)\omega'}\Re\left(\sum_{n<X^2}\frac{\Lambda_X(n)}{n^{\sigma_1+ti}}\right)\nonumber \\
&~~~~~~~+\frac{\left(1+\frac{1}{e}\right)\frac{1}{e}\omega}{1-\frac{1}{e}\left(1+\frac{1}{e}\right)\omega'}\cdot\frac{1}{2}\log t+O\left(X^{\frac{1}{2}-\sigma_1}\right)\Biggr\}\nonumber\\
&=\frac{1}{(m+1)(\log X)^{m+1}}\cdot \frac{\left(1+\frac{1}{e}\right)\frac{1}{e}\omega}{1-\frac{1}{e}\left(1+\frac{1}{e}\right)\omega'}\cdot \frac{1}{2} \log t \nonumber \\
&~~~~~~~~~~~~~~~~~~~~~~+O\left\{\frac{1}{(\log X)^{m+1}}\left|\sum_{n<X^2}\frac{\Lambda_X(n)}{n^{\sigma_1+ti}}\right|\right\}\nonumber \\
&\ll\eta'_3(t)+\frac{X}{(\log X)^{m+2}},\label{1.17}
\end{align}
say. As well as $\eta_1(t)$, we have
\begin{align*}
|\eta'_1(t)|\leq \frac{1}{1-\frac{1}{e}\left(1+\frac{1}{e}\right)}\cdot \frac{1}{2}\log t \cdot \frac{1}{(\log X)^{m+1}}\left(\sum_{j=0}^{m}\frac{m!}{(m-j)!}\left(\frac{1}{e}+\frac{1}{2^{j+1}e^2}\right)\right).
\end{align*}
Finally, we estimate $J_3$. By Stirling's formula, we get
\begin{align}
\left|\frac{\Gamma'}{\Gamma}\left(\frac{\sigma_1+ti}{2}+1\right)\right|&=\left|\frac{i}{2}\log\frac{ti}{2}+\left(\frac{\sigma_1+ti+1}{2}\right)\frac{1}{t}-\frac{i}{2}+O\left(\frac{1}{t}\right)\right|\nonumber \\
&\leq\frac{1}{2}\log t+O\left(\frac{1}{t}\right).\label{1.18}
\end{align}
Also $\left|\frac{\Gamma'}{\Gamma}\left(\frac{\sigma+ti}{2}+1\right)\right|$ is estimated similarly.

Hence by (\ref{1.18}) and Lemma \ref{lem2}, we have
\begin{align*}
&\left|\Im\left\{\frac{\zeta'}{\zeta} (\sigma_1+ti )-\frac{\zeta'}{\zeta} (\sigma+ti )\right\} \right | ~~~~~~~~~~~~~~~~\\
&~~~~~~~~~~~~~~\leq \Im\left\{\sum_{\rho}\left(\frac{1}{\sigma_1+ti-\rho}+\frac{1}{\rho}\right)-\sum_{\rho}\left(\frac{1}{\sigma+ti-\rho}+\frac{1}{\rho}\right)\right\}+O\left(\frac{1}{t}\right)\\
&~~~~~~~~~~~~~~=\sum_{\gamma}\frac{(t-\gamma)\left\{\left(\sigma-\frac{1}{2}\right)^2-\left(\sigma_1-\frac{1}{2}\right)^2\right\}}{\left\{\left(\sigma_1-\frac{1}{2}\right)^2+(t-\gamma)^2\right\}\left\{\left(\sigma-\frac{1}{2}\right)^2+(t-\gamma)^2\right\}}+O\left(\frac{1}{t}\right).
\end{align*}
Therefore,
\begin{align*}
|\Im(J_3)|&\leq \left|\int_{\frac{1}{2}}^{\sigma_1}\left(\sigma-\frac{1}{2}\right)^m \sum_{\gamma}\frac{(t-\gamma)\left\{\left(\sigma-\frac{1}{2}\right)^2-\left(\sigma_1-\frac{1}{2}\right)^2\right\}}{\left\{\left(\sigma_1-\frac{1}{2}\right)^2+(t-\gamma)^2\right\}\left\{\left(\sigma-\frac{1}{2}\right)^2+(t-\gamma)^2\right\}}d\sigma\right|\\
&~~~+\int_{\frac{1}{2}}^{\sigma_1}\left(\sigma-\frac{1}{2}\right)^m\cdot O\left(\frac{1}{t}\right)d\sigma.
\end{align*}
If $t=\gamma$, the first term of the right-hand side of above inequality is $0$. If $t\neq \gamma$, since $\sigma<\sigma_1$, we have
\begin{align*}
&\left|\int_{\frac{1}{2}}^{\sigma_1}\left(\sigma-\frac{1}{2}\right)^m \left\{\sum_{\gamma}\frac{(t-\gamma)\left\{\left(\sigma-\frac{1}{2}\right)^2-\left(\sigma_1-\frac{1}{2}\right)^2\right\}}{\left\{\left(\sigma_1-\frac{1}{2}\right)^2+(t-\gamma)^2\right\}\left\{\left(\sigma-\frac{1}{2}\right)^2+(t-\gamma)^2\right\}}\right\}d\sigma \right| \\
&~~~<\sum_{\gamma}\frac{\left(\sigma_1-\frac{1}{2}\right)^{m+2}}{\left(\sigma_1-\frac{1}{2}\right)^2+(t-\gamma)^2} \int_{\frac{1}{2}}^{\sigma_1}\frac{|t-\gamma |}{\left(\sigma-\frac{1}{2}\right)^2+(t-\gamma)^2}d\sigma \\
&~~~\leq \sum_{\gamma}\frac{\left(\sigma_1-\frac{1}{2}\right)^{m+2}}{\left(\sigma_1-\frac{1}{2}\right)^2+(t-\gamma)^2} \int_{\frac{1}{2}}^{\infty}\frac{|t-\gamma|}{\left(\sigma-\frac{1}{2}\right)^2+(t-\gamma)^2}d\sigma \\
&~~~\leq \frac{\pi}{2}\left(\sigma_1-\frac{1}{2}\right)^{m+1}\sum_{\gamma}\frac{\sigma_1-\frac{1}{2}}{\left(\sigma_1-\frac{1}{2}\right)^2+(t-\gamma)^2}.
\end{align*}
Applying (\ref{1.9}) and (\ref{sigxx}), and taking $X=\log t$ lastly, the right-hand side of above inequality is
\begin{align}
&~~~\leq \frac{\pi}{2}\left(\sigma_1-\frac{1}{2}\right)^{m+1}\left\{\frac{1}{1-\frac{1}{e}\left(1+\frac{1}{e}\right)\omega'}\cdot\frac{1}{2}\log t+O\left(\left|\sum_{n<X^2}\frac{\Lambda_X(n)}{n^{\sigma_1+ti}}\right|\right)\right\}\nonumber\\
&~~~=\frac{\pi}{4}\cdot \frac{1}{1-\frac{1}{e}\left(1+\frac{1}{e}\right)\omega'}\cdot \frac{\log t}{(\log X)^{m+1}}+O\left(\frac{X}{(\log X)^{m+2}}\right)\nonumber\\
&~~~\leq \frac{\pi}{4}\cdot \frac{1}{1-\frac{1}{e}\left(1+\frac{1}{e}\right)}\cdot \frac{\log t}{(\log\log t)^{m+1}}+O\left(\frac{\log t}{(\log\log t)^{m+2}}\right).\label{1.20}
\end{align}
Also,
\begin{align}
\int_{\frac{1}{2}}^{\sigma_1}\left(\sigma-\frac{1}{2}\right)^m\cdot O\left(\frac{1}{t}\right)d\sigma=O\left(\frac{1}{t(\log X)^{m+1}}\right). \label{1.21}
\end{align}
By (\ref{1.20}) and (\ref{1.21}),
\begin{align}
|\Im(J_3)| &\leq \frac{\pi}{4}\cdot \frac{1}{1-\frac{1}{e}\left(1+\frac{1}{e}\right)}\cdot \frac{\log t}{(\log \log t)^{m+1}}\nonumber \\
&~~~~~~~~~~~~~~~~~~~~~~~~
+O\left(\frac{1}{t(\log\log t)^{m+1}}\right)+O\left(\frac{\log t}{(\log\log t)^{m+2}}\right).\label{F}
\end{align}
Therefore, we obtain by (\ref{1.15}), (\ref{1.16}), (\ref{1.17}), (\ref{F}), $\eta'_1(t)$ and $\eta'_3(t)$
\begin{align*}
|S_m(t)|&\leq \frac{1}{2\pi m!}\cdot\frac{\log t}{(\log \log t)^{m+1}}\Biggl\{\frac{1}{1-\frac{1}{e}\left(1+\frac{1}{e}\right)}\sum_{j=0}^{m}\frac{m!}{(m-j)!}\left(\frac{1}{e}+\frac{1}{2^{j+1}e^2}\right)\\
&~~+
\frac{1}{m+1}\cdot \frac{\left(1+\frac{1}{e}\right)\frac{1}{e}}{1-\frac{1}{e}\left(1+\frac{1}{e}\right)} +
\frac{\pi}{2}\cdot \frac{1}{1-\frac{1}{e}\left(1+\frac{1}{e}\right)}\Biggr\}+O\left(\frac{\log t}{(\log\log t)^{m+2}}\right).
\end{align*} 
\qed
\\
{\bf Acknowledgments}\\

I thank my advisor Prof. Kohji Matsumoto for his advice and patience during the preparation of this paper. I also thank members in the same study, who gave adequate answers to my questions. Finally, I thank referees who indicate errors in this paper.

\end{document}